\documentclass[a4paper,11pt]{amsart}
\usepackage{t1enc,amssymb,amsmath}
\usepackage[english]{babel}

\input xy
\xyoption{all}

\newtheorem{theorem}{Theorem}[section]
\newtheorem{corollary}[theorem]{Corollary}

\newtheorem{proposition}[theorem]{Proposition}

\theoremstyle{remark}
\newtheorem{remark}[theorem]{\sc Remark}
\theoremstyle{definition}
\newtheorem{definition}[theorem]{Definition}

\theoremstyle{remark}

\renewcommand{\Box}{_\square}    

\newcommand{\gras}[1]{{\mathbb #1}}
\newcommand{\C}{\gras{C}}
\newcommand{\Z}{\gras{Z}}

\newcommand{\R}{\gras{R}}
\newcommand{\T}{\gras{T}}

\begin{document}

\title[Fillability and cohomology]{\textbf{On the cohomology rings \\ 
  of holomorphically fillable manifolds}}

\author{Patrick Popescu-Pampu}
\address{Univ. Paris 7 Denis Diderot, Inst. de
  Maths.-UMR CNRS 7586, {\'e}quipe "G{\'e}om{\'e}trie et dynamique" \\ 
  Site Chevaleret, Case
  7012, 75205 Paris Cedex 13, France.}
\email{ppopescu@math.jussieu.fr}

\date{19 December 2007}


\thispagestyle{empty}
\begin{abstract}{An odd-dimensional differentiable manifold is called 
   \emph{holomorphically fillable} if it is diffeomorphic to the boundary of 
   a compact strongly pseudoconvex complex manifold, \emph{Stein
     fillable} if this last 
   manifold may be chosen to be Stein and \emph{Milnor fillable} if it is 
   diffeomorphic to the abstract boundary of an isolated singularity of 
   normal complex analytic space. We show that the homotopical 
   dimension of a  manifold-with-boundary of dimension at least $4$  
    restricts the cohomology ring (with any coefficients) of its boundary. 
    This gives restrictions on the cohomology rings of Stein fillable 
    manifolds, on the dimension of the exceptional locus of any  resolution 
    of a given isolated  singularity, and on the topology of smoothable 
    singularities.  We give also new proofs of structure theorems of Durfee \&
    Hain and Bungart about the cohomology rings of Milnor fillable and
    respectively holomorphically fillable manifolds. The various
    structure theorems presented in this paper imply that in dimension
  at least 5, the classes of Stein fillable, Milnor fillable and
  holomorphically fillable manifolds are pairwise different.}   
\end{abstract}

\maketitle

\par\medskip\centerline{\rule{2cm}{0.2mm}}\medskip
\setcounter{section}{0}

\section{Introduction}

The foundational papers \cite{G 85}, \cite{E 89}, \cite{E 90}, \cite{EG 91} 
of Eliashberg and Gromov showed   that one can 
get information on the structure of a contact manifold $N$, 
whenever this manifold bounds an even dimensional manifold 
$W$ with a holomorphic or symplectic structure compatible in some way 
with the contact structure on the boundary: one says that $N$ is 
\emph{filled} by $W$. Since then, many notions of \emph{fillability} 
for contact manifolds have been introduced: \emph{holomorphic, Stein,
  Milnor, 
Liouville, Weinstein, strong symplectic, weak symplectic, etc} (see
Geiges \cite{G 08}). 

In this paper we restrict to the notions of holomorphic, Stein and 
Milnor fillability. An odd-dimensional, closed, orientable  manifold 
is called \emph{holomorphically fillable} if it is 
diffeomorphic to the boundary of a compact 
strongly pseudoconvex complex manifold. It is called \emph{Stein
  fillable} if this manifold can  
be chosen to be Stein. It is called \emph{Milnor fillable} if it is
diffeomorphic to the abstract boundary (or link) of an isolated
singular point of normal complex analytic space. Using a resolution of
the singularity, we see that a Milnor fillable manifold is
automatically holomorphically fillable. Notice that we have
dropped the contact 
structures from these definitions, due to the fact that in this article 
we prove theorems which involve only  
the cohomology rings of each kind of fillable manifolds. 
Nevertheless, we 
began this article by speaking about them because it is the 
interest in Milnor fillable \emph{contact} manifolds  
which led us to those theorems, before discovering that some of them
were already known. 

Holomorphically fillable 3-manifolds are necessarily 
Stein fillable, as was proved by Bogomolov \& de Oliveira \cite{BO 97}.  
In higher dimensions, this is no longer the case. For example,
Eliashberg, Kim \& Polterovich \cite{EKP 06} have shown that the
projective spaces $\R\gras{P}^{2n-1}$, which are always
holomorphically fillable, are not Stein fillable whenever $n \geq 3$.

\medskip
In Section \ref{notions} we give the background on (strictly)
plurisubharmonic functions, strongly pseudoconvex spaces, Stein
spaces, fillable manifolds and resolutions of isolated singularities
of normal complex analytic spaces needed in the rest of the paper. 

In  Section \ref{constr}, we 
show that the homotopical dimension of a compact
manifold-with-boundary  restricts the 
cohomology ring with arbitrary coefficients of its boundary (Theorem
\ref{hyphom}).  This extends the method used by  
Eliashberg, Kim \& Polterovich to show that $\R\gras{P}^{2n-1}$ is not
Stein fillable. A Stein  manifold being homotopically 
of dimension at most equal to its complex dimension, we get in 
particular constraints on the cohomology rings of Stein fillable 
manifolds (Corollary \ref{Stein}). Then we consider manifolds 
$N$ which are total spaces of oriented circle bundles $N 
\stackrel{p}{\longrightarrow} \Sigma$, and we apply 
Theorem \ref{hyphom} by showing  that suitable 
hypotheses on the cohomology ring of  $\Sigma$ and on the Euler class 
of the bundle, give lower bounds on the homotopical dimension 
of any filling of $N$ (Proposition \ref{bundles}).  

In  Section \ref{appl}, we give applications of the results of the 
previous section to isolated singularities of complex analytic spaces. 
First, we give a lower bound on the dimension of the exceptional 
locus of any resolution in terms of the  cohomology ring of the 
boundary (Proposition \ref{dimex}). 
The generic fibers of a smoothing of an 
isolated singularity being Stein and their boundaries being  
diffeomorphic to the boundary of the singularity, we also get 
constraints on the topology of smoothable singularities 
(Proposition \ref{nonsm}).  We consider 
in more detail the isolated singularities  
obtained by contracting the zero-section of an  
anti-ample line bundle $L$  on a projective manifold $\Sigma$.  Suitable
hypotheses on the integral cohomology ring $H^*(\Sigma; \Z)$ and on
the Chern class of $L$ imply that the boundary of the resulting  singularity
(Milnor fillable, therefore holomorphically fillable) is not Stein fillable.  
In particular, such a germ is \emph{non-smoothable}. 
As a special case  we get (Corollary \ref{corfin}):

\medskip
 \emph{Let $(X,x)$ be the germ of normal analytic space with isolated
   singularity  
 obtained by contracting the $0$-section of an anti-ample line bundle on 
 an abelian variety $\Sigma$ of complex dimension $\geq 2$, and whose 
 first Chern class is not primitive in $H^2(\Sigma, \Z)$. Then the 
 boundary of $(X,x)$ is  
 not Stein fillable. In particular, $(X,x)$ is not smoothable. }
 \medskip
 
 We answer like this partially the concluding question asked by Biran in 
 \cite{B 05} (see Remark \ref{biran}).

In Section \ref{Milnor} we give a new proof of a theorem of Durfee \&
Hain \cite{DH 88} (first announced in \cite{D 86}), describing
restrictions on the cohomology rings with
rational coefficients of Milnor fillable manifolds of arbitrary
dimension (Theorem \ref{obsmiln}). This  generalizes a
theorem obtained by Sullivan \cite{S 75} in dimension 3. A crucial ingredient 
in our proof is a theorem of Goresky \& MacPherson \cite{GM 82},  
describing the kernel of the 
map from the homology of the boundary to the homology of the manifold, in 
the case of a divisorial resolution of an isolated singularity 
(Theorem \ref{decomp}).  This theorem generalizes to any dimension the 
fact that the intersection form associated to a resolution of a normal 
surface singularity is non-degenerate. Both our proof and the one of
Durfee \& Hain are based on a deep purity theorem of Beilinson,
Bernstein, Deligne \& Gabber \cite{BBD 82} (in our proof, this is
hidden inside Goresky \& MacPherson's theorem).

 In Section \ref{holomorph} we  apply the results of Section
 \ref{Milnor}  in order to give a new proof of a theorem of Bungart
 \cite{B 92}, showing that the rational cohomology rings of
holomorphically  fillable manifolds are also constrained (Theorem
\ref{holofill}).  
 
 We deduce examples in any odd dimension 
$\geq 3$ of holomorphically fillable manifolds which are not Milnor fillable 
(Corollary \ref{prod}). Combining them with the examples of the previous 
sections, we see that in all odd dimensions $\geq 5$, the classes of 
Stein, Milnor or holomorphically fillable manifolds are pairwise distinct.

\medskip
\begin{small}
\textbf{Acknowledgments.} I am grateful to Yakov Eliashberg, who
gave me the reference  
\cite{EKP 06} after receiving a previous version of this
paper; this made me simplify the criterion of non-fillability by
manifolds having small homotopical dimension. He informed me that  
he had learnt from Freedman  at the beginning of the 1990's 
the possibility to get obstructions on the cohomology rings 
of the boundaries of manifolds which are homotopically of dimension
equal to half  their real dimension. I had the same idea inspired by 
the work \cite{S 75} of Sullivan. I am grateful to Etienne Ghys who 
showed me Sullivan's article a few years ago. This made me find the 
results of this paper before discovering that some of them had already 
been proved by Durfee \& Hain and Bungart. 
I am also grateful to Hansj{\"o}rg Geiges, Eduard Looijenga, David
Mart{\'\i}nez, Luca Migliorini, Jan Schepers, Jos{\'e} Seade and  
Bernard Teissier for stimulating conversations. 
\end{small}

\medskip
\section{Plurisubharmonic functions and 
various notions of  fillability} \label{notions} 

For details on the various notions and results recalled in this
section, one can consult  
Grauert \& Remmert \cite{GR 79}, Peternell \cite{P 94}, Bennequin
\cite{B 90} and  Eliashberg \cite{E 97}.

Let $X$ be a complex manifold. Denote by $TX$ the (real) tangent bundle of the 
underlying real differentiable manifold and by $J:TX\rightarrow TX$ the 
(integrable) almost complex structure associated to the complex
structure of $X$.  
The operator $d^{\C}:=J^*\circ d$ (that is, $d^{\C}f:= df \circ J$ for
any smooth  
function defined on $X$) is real and  intrinsically associated to the 
complex structure of $X$. In terms of the operators $\partial=d'$ and 
$\overline{\partial}=d''$, one has:
$$\left\{ \begin{array}{l}
                   d=\partial + \overline{\partial}\\
                   d^{\C}=i(\partial - \overline{\partial})
                \end{array} \right. .$$  

Let $f$ be a real-valued smooth function defined 
on $X$.  

We associate to $f$ the following tensors on $X$:
$$\begin{array}{l}
        \alpha_f:=-d^{\C}f,\\
        \omega_f:= d\alpha_f= -dd^{\C}f,\\
        g_f(u,v):=\omega_f(u, Jv), \: \forall \: u,v \in TV,\\
        h_f :=g_f +i \omega_f.
      \end{array}$$ 
 The kernel of the restriction of $\alpha_f$ to any regular level 
 $X_a:=f^{-1}(a)$ of $f$ is the 
 complex tangent bundle $TX_a\cap J(TX_a)$ of $X_a$. 
The exterior form $\omega_f$ is  real of type $(1,1)$. It is the  associated 
\emph{Levi form} of the function $f$. The associated hermitian form $h_f$ 
($\C$-antilinear in the first coordinate and $\C$-linear in the second
coordinate)  
is also called the \emph{Levi form} of $f$. If needed, we distinguish 
between the two versions of Levi form by precising that we deal with 
\emph{the exterior Levi form} or \emph{the hermitian Levi form}.  

In the sequel we will be mainly interested in a special class of real-valued 
functions on $X$:

\begin{definition} \label{psh}
The function $f:X \rightarrow \R$ is called 
\textbf{plurisubharmonic} (abbreviated \textbf{psh}) if $g_f$ is positive 
semidefinite. It is called \textbf{strictly plurisubharmonic} 
(abbreviated \textbf{spsh}) if $g_f$ is positive definite, that is, if it 
defines a riemannian structure on the smooth manifold $X$. 
\end{definition}

Notice that $f$  is spsh if and only if the associated Levi form is
K{\"a}hler. (Strictly)  
plurisubharmonic functions are the analogs on complex manifolds of 
(strictly) convex functions on real manifolds endowed with an affine 
structure (not to be confused with affine algebraic manifolds). 

The notion 
of (s)psh function can be defined also  if $X$ is a reduced but not
necessarily  smooth complex space: 
the  function $f:X \rightarrow \R$ is called \emph{(s)psh} if in a
neighborhood of  
each point of $X$, it is the restriction to $X$ of a smooth (s)psh function 
defined on a complex manifold into which one has locally embedded $X$. This 
definition does not depend on the choice of local embedding.

 Strictly plurisubharmonic functions have the following easily provable 
 properties:

\begin{proposition} \label{spsh}

1) The restriction of a (s)psh function to a complex subspace is again (s)psh.

2) If $f$ is (s)psh and $\phi: \R \rightarrow \R$ is (strictly) convex
and smooth,     then  $\phi\circ f: X \rightarrow \R$ is (s)psh. 
     
3) spsh functions form an open set among smooth functions 
 in the $C^2$ topology with compact supports.
\end{proposition}

The notions of (s)psh functions are local. The following notion is instead 
global: the smooth real-valued function $f$ defined on the reduced complex 
analytic space $X$ is called an \emph{exhaustion function}  if it is
proper and  bounded from below (which is equivalent to the fact that
it attains its  
absolute minimum). 

\begin{definition}
The reduced complex analytic space $X$ is called \textbf{strongly 
pseudoconvex} if it carries an exhaustion function $f$  which is strictly 
pluri\-subharmonic outside a compact set. $X$ is called \textbf{a Stein 
space} if $f$ may be chosen to be strongly pseudoconvex all over $X$.
\end{definition} 

The following  characterizations of Stein spaces may be obtained by combining  
theorems of 
Grauert and Narasimhan (see Grauert \& Remmert \cite[page 152]{GR 79} or 
Peternell \cite[sections 1--4]{P 94}):

\begin{theorem} \label{propstein} 
Let $X$ be a reduced paracompact complex
  analytic space. The following properties are
  equivalent: 

1) $X$ is a Stein space.

2) $X$ is holomorphically convex and holomorphically separable (that
is, the global holomorphic functions separate the points). 

3) $X$ is strongly pseudoconvex and has no compact 
analytic subsets of positive dimension.
\end{theorem}

There are also characterizations of Stein spaces using coherent cohomology, 
but we won't use them in this paper.  As a corollary of the previous theorem, 
we get:

\begin{theorem} \label{fin}
  If $X\rightarrow Y$ is a finite map and $Y$ is a 
  Stein space,    then $X$ is also Stein.
\end{theorem}

In particular, any closed subspace of a Stein space is Stein. Still more 
particularly, any closed subspace of $\C^n$ is Stein. In fact, in this way 
one does not restrict very much the class of Stein spaces. Indeed, as a 
result of works of Remmert, Bishop and Narasimhan, one has the 
following embedding theorem (see Bell \& Narasimhan 
\cite[Theorem 3.1]{BN 90}):

\begin{theorem} \label{embed}
  A Stein space $X$ can be embedded holomorphically in some $\C^n$ 
  if and only if it has bounded local embedding dimension (that is, the 
  dimension of the Zariski tangent spaces is a bounded function on $X$).
\end{theorem}

Suppose that $X$ is strongly pseudoconvex. Then it has a maximal 
compact analytic subset $K\subset X$ with a finite number of 
irreducible components. Consider the Remmert reduction morphism 
$X \stackrel{r}{\rightarrow} Y$ (see Peternell 
\cite[Section 2]{P 94}). It contracts to a point each \emph{connected} 
component of $K$ and it is an isomorphism outside $K$.  
The space $Y$ is then a Stein space, by point 1) of the previous theorem. 
If $X$ is a manifold, then $Y$ is a normal Stein space with a finite number 
of isolated singularities. For this reason,  
in the sequel we will consider only complex spaces with isolated
singularities.  
If $f$ is an exhaustion function  on the complex space $X$ and
$X_a:=f^{-1}(a)$ is  
a regular level of $f$ in whose neighborhood $f$ is spsh, then we say that 
the compact sublevel $X_{\leq a}:= f^{-1}((-\infty, a])$ is \emph{a
  compact strongly  
pseudoconvex manifold}, whose boundary is $X_a$. If $X$ is Stein, we say that 
$X_{\leq a}$ is \emph{a compact Stein manifold}.

We will do Morse theory starting from  spsh functions 
defined on complex manifolds $X$. The fundamental observation  is:

\begin{proposition} \label{index}
 If the complex manifold $X$ has complex dimension $n \geq 1$, then the 
 indices of  the critical  points of  a spsh Morse function on $X$ are
 $\leq n$.  
\end{proposition}

This remark made by Thom around 1957 was the starting point of proofs 
through Morse theory of Lefschetz' hyperplane section theorem, by 
Andreotti \& Frankel \cite{AF 59} and Bott \cite{B 59}. Milnor \cite{M
  63} noticed that an immediate consequence of those proofs is:

\begin{theorem} \label{homtype}
  A Stein manifold $X$ of complex dimension $n$ has the homotopy type of a 
  CW-complex of   dimension at most $n$. As a consequence, all the homology 
  and cohomology groups with arbitrary 
 coefficients of $X$ vanish in degree at least $n+1$. 
\end{theorem}

The vanishing of the previous cohomology groups with real coefficients 
was first proved by Serre \cite{S 53}. An analogous vanishing theorem 
was proved for arbitrary Stein spaces by Narasimhan \cite{N 67}. 
The analog of the first sentence 
of the theorem was then proved for affine algebraic spaces by Karchyauskas 
\cite{K 79} and for arbitrary Stein spaces by Hamm \cite{H 86}. In
this paper we will need only the (co)homological version of the
theorem for Stein manifolds, but (and this is very important)
\emph{with integral coefficients}. We would like also to mention that in 
\cite{E 90}, Eliashberg characterized the differentiable manifolds of 
even dimension $\geq 6$ which 
admit a Stein structure. 

In \cite{E 97}, 
Eliashberg explains Proposition \ref{index} using symplectic geometry in 
the following way, which is excellent  for understanding
the interdependence of the objects $\alpha_f, \omega_f, g_f$ defined
before. Consider the gradient of $f$ with respect to the riemannian  
metric $g_f$. It is also equal to the Liouville vector field of
$\alpha_f$ with respect  
to the symplectic form  $\omega_f$, therefore its flow exponentially
dilates this  symplectic form. This implies that the stable cells of
the gradient 
associated to  the critical points of $f$ are isotropic (that is,
$\omega_f$ vanishes 
on them),  therefore their dimension is at most $n$.

In this article, the most important example of spsh function is 
the squared-distance function to a point in some space $\C^n$ (which
was also the type of function used by Andreotti \& Frankel \cite{AF
  59}). We will also  
consider restrictions of such functions to complex analytic subspaces
of $\C^n$, which are again spsh by Proposition \ref{spsh}.

In particular, let $(X,x)$ be an irreducible  germ of reduced complex
analytic space  
which is smooth outside $x$. We also say that $(X,x)$ is \emph{an
  isolated singularity}.  
Choose an embedding of a representative of $(X,x)$ in some $(\C^n,0)$. 
Denote by $\rho:(X,x) \rightarrow (\R_+,0)$ the restriction of the
squared-distance  
function to the origin. We call it \emph{a euclidean rug function}
associated to the  
isolated singularity $(X,x)$ (notice that in \cite{CNP 06} we had
introduced more general  
euclidean rug functions; the name is inspired by Thom's article \cite{T 59}). 
Its levels $\rho^{-1}(\epsilon)$ are all smooth and diffeomorphic for
$\epsilon \in (0, \epsilon_0]$, where $\epsilon_0>0$ is sufficiently
small. Their diffeomorphism type  
does not depend on the choice of the embedding, and is called 
\emph{the (abstract) boundary} or \emph{the (abstract) link} of
$(X,x)$. We say that a compact representative of $(X,x)$ of the form
$\rho^{-1}([0, \epsilon_0])$ with the properties stated before is
\emph{a compact Milnor representative of} $(X,x)$. This notion may be
extended to any reduced germ, without necessarily isolated
singularity.

\begin{definition}
  The odd-dimensional manifold $N$ is called \textbf{holomorphically fillable} 
  if it is diffeomorphic to the boundary of a compact strongly
  pseudoconvex space $X$  with at most 
  isolated singularities. $N$ is called 
  \textbf{Stein fillable} if $X$ may be chosen to be a Stein \emph{manifold}. 
  $N$ is called \textbf{Milnor fillable} if it is diffeomorphic to
  the abstract boundary  of an isolated singularity. 
\end{definition} 

Note that by Hironaka's theorem of resolution of singularities and the
fact that  the Remmert reduction of a strongly pseudoconvex space is a
bimeromoprhism, we would 
obtain equivalent  definitions of holomorphic fillability by asking 
$X$ to be either smooth or to be Stein with at most isolated singularities. 

The notions of holomorphic fillability and Stein fillability were
introduced in the context of the study of convexity notions in
symplectic geometry 
by Eliashberg $\&$ Gromov \cite{EG 91}. The notion of
Milnor fillability was introduced by Caubel, N{\'e}methi and myself 
in the paper \cite{CNP 06}. In all these cases, one considers a
supplementary \emph{contact} structure on the manifold $N$ (see also
Geiges \cite{G 08}), and one has to take care of orientation
issues. As in this article we give purely cohomological obstructions
for Stein, Milnor and holomorphic fillability, we do not spend time
here on these issues. We precise only, for the reader who wants to get
an idea of the relation between what we explained before and contact
geometry, that whenever $f$ is spsh, the restriction of the real
$1$-form $\alpha_f$ to a regular level $f^{-1}(a)$ of $f$ is a contact
form, and 
that the orientation defined on this level by $\alpha_f \wedge
(d\alpha_f)^{\wedge (n-1)}$ coincides with its orientation as a
boundary of the sublevel $f^{-1}((-\infty, a])$. Varchenko \cite{V 80} 
showed that the associated contact structures on the boundaries of 
Milnor representatives of isolated singularities are independent of 
the choice of euclidean rug function. In \cite{CNP 06}, we continued the 
study of such \emph{Milnor fillable contact manifolds}. 

In Section \ref{holomorph} we will study holomorphically fillable
manifolds by considering a 
convenient cobordism which relates them to a disjoint union of Milnor fillable
manifolds. This cobordism will be constructed using the
following proposition:

\begin{proposition} \label{specfunc}
   Let $X$ be a compact Stein space with isolated singularities and 
   $f:X\rightarrow \R$ a spsh exhaustion function. Denote by $F$ a finite set 
   which contains the singular locus of $X$. Then there exists a spsh function 
   $\phi :X \rightarrow \R$ which coincides with $f$ outside a compact
   subset of  
   the interior of $X$, which attains its absolute minimum exactly on $F$ and 
   which is a euclidean rug function in restriction to a sufficiently small 
   neighborhood of any point of $F$ in $X$. 
\end{proposition} 

\textbf{Proof.} There exists an analytic morphism $X \rightarrow Y$ which
identifies  all the points of $F$ and is a biholomorphism outside $F$. 
Denote by $y \in Y$ the image of $F$ by this morphism. 
Point 3) of Theorem \ref{propstein} implies that  $Y$ is again a
compact Stein space.  
By post-composing $f$ with an adequate  smooth convex function which 
interpolates between a constant function and the identity on $\R$, we get 
a psh function $\tilde{f}$ on X which is constant on $F$ and equal to $f$ 
near $\partial X$. Therefore it descends to a psh function on $Y$.

As the local embedding dimension on $Y$ is globally bounded, by
Theorem \ref{embed} the interior 
$\stackrel{\circ}{Y}$ of $Y$ can be embedded in some space $\C^n$.  
From now on,  we shall look at 
$Y$ as a subspace of $\C^n$.  Consider then the function $\rho:Y
\rightarrow \R$  
obtained by restricting to $Y$ the squared-distance function to $y$ in 
$\C^n$. 

Consider a number $a >0$ such that $\tilde{f}=f$ wherever $\rho >a$
(that is, such  
that we have modified $f$ only inside the euclidean ball with radius
$\sqrt{a}$  
centered at $y$). Then post-compose $\rho$ with a smooth real-valued 
function defined on $\R$, which is the identity on the interval
$]-\infty, a]$ and  
identically zero on $[b, +\infty[$, where $b>a$ is chosen
arbitrarily. Notice that we  
impose nothing more than its global smoothness in between. Denote by 
$\tilde{\rho}:Y \rightarrow \R$ the function obtained like this. It is
a euclidean rug  
function in a neighborhood of $y$, it is spsh wherever $\tilde f \neq f$ 
and it vanishes outside a compact. 

By using Proposition \ref{spsh}, we see that the function $\psi_{\epsilon}:=
\tilde{f} + \epsilon \cdot \tilde{\rho}$ is spsh all over $Y$ whenever
$\epsilon>0$  
is sufficiently small. Then its lift $\phi_{\epsilon}$ to $X$
verifies all the  
conditions asked for in the conclusion of the proposition. \hfill $\Box$

\medskip
As a preliminary to the study of holomorphically fillable manifolds,
in Section  \ref{Milnor} we will prove a theorem of Durfee \& Hain,
showing  that there exist restrictions on
the rational cohomology rings of Milnor fillable manifolds of
dimension at least 3. Our proof is different form the one of Durfee \&
Hain. It is based on a theorem of
Goresky \& MacPherson. Let us explain it. 

Consider an isolated singularity $(X,x)$ of normal complex analytic
  space. A \emph{resolution} of $(X,x)$ is a proper morphism  $(\tilde{X},
E)\longrightarrow (X,x)$ with smooth total space $\tilde{X}$,  
   realizing an isomorphism outside the singular locus $x$. By
   Hironaka's theorem, resolutions exist. A resolution is called
   \emph{divisorial} if its exceptional set $E$ is purely of codimension
   1 in $\tilde{X}$. All the resolutions of normal surfaces are
   divisorial, but this is not true in higher dimensions (the simplest
   example of non-divisorial resolution is recalled at the beginning
   of Section \ref{appl}).  

  Suppose now that $X$ denotes a compact Milnor representative of
  the germ. Let  $(W,E)\stackrel{\pi}{\rightarrow}(X,x)$ be a divisorial
  resolution whose exceptional divisor $E$ has only normal
  crossings. Denote by $N$ 
  the boundary of $W$. As $\pi$ is an isomorphism outside $x$, it 
  identifies $N$ with the boundary of
  $X$, that is with the abstract boundary of the singularity.  
  We will need the following theorem relating the (co)homology of
  the boundary $N$ of the singularity to the (co)homology of the
  resolution $W$:

\begin{theorem}\label{decomp}
The following are equivalent facts and are true if $W$ is a divisorial
resolution of a Milnor representative of a normal isolated singularity
and $N$ is its boundary (homology and cohomology groups  are
considered with \emph{rational} coefficients and all the morphisms are
induced by inclusions): 

1) The morphisms $H_i(N)\rightarrow H_i(W)$ vanish identically for $i
\in \{n,...,\linebreak 2n-1\}$. 

2) The morphisms $H_i(N)\rightarrow H_i(W)$ are injective for $i
\in \{1,...,n-1\}$ and vanish identically for $i
\in \{n,...,2n-1\}$. 

3) The morphisms $H^i(W)\rightarrow H^i(N)$ vanish identically for $i
\in \{n,...,\linebreak 2n-1\}$. 

4) The morphisms $H^i(W)\rightarrow H^i(N)$ are surjective for $i \in
\{0,...,n-1\}$ and vanish identically for $i
\in \{n,...,2n-1\}$. 

5) The morphisms $H_i(W)\rightarrow H_i(W,N)$ are surjective for $i
\in \{1,...,n\}$ and are injective for $i \in \{n,..., 2n-1\}$. 

6) The morphisms $H^i(W,N) \rightarrow H^i(W)$ are
injective for $i \in \{1,...,n\}$ and surjective for $i
\in \{n,..., 2n-1\}$. 
\end{theorem}

The equivalence between 1)--6) results from a play with the long
exact (co)homology sequences of the pair $(W,N)$, Poincar{\'e}-Lefschetz
duality and the fact that $H^i(W)\rightarrow H^i(N)$ is the adjoint
morphism of $H_i(N)\rightarrow H_i(W)$ when we work over
$\mathbf{Q}$. 
Point 1) of the theorem was deduced by Goresky \& MacPherson in
\cite[page 124]{GM 82} as a 
consequence of a deep decomposition theorem in intersection homology
theory proved by Beilinson, Bernstein, Deligne \& Gabber
\cite[Theorem 6.2.5, page 163]{BBD 82}. Point 6) 
was proved by Steenbrink
\cite[page 518]{S 83} as part of a short exact sequence of mixed Hodge
structures (see also \cite[page 117]{S 87}). A
Hodge-theoretic proof of this theorem was given by Navarro Aznar in
\cite[page 285]{N 85}. It can also be obtained as a consequence of de
Cataldo $\&$ Migliorini's \cite[Corollary 2.1.12]{CM 05}. 

As noted in \cite[page 123]{GM 82}, for any compact oriented manifold
$W$ with boundary $N$,  the kernel of the morphism 
$H_{\bullet}(N)\rightarrow H_{\bullet}(W)$ between the total 
homologies, induced by the inclusion $ N \hookrightarrow W$, is
half-dimensional inside $H_{\bullet}(N)$. The previous theorem
describes this kernel when $W$ is a divisorial resolution of an
isolated singularity: it is exactly $\oplus_{i=n}^{2n-1} H_i(N)$. In
fact, this was the way in which Goresky \& MacPherson stated their
theorem. 

In the case of a germ of surface $(X,x)$, the previous theorem is a
consequence of the fact (proved by Du Val \cite{V 44} and
Mumford \cite{M 61}) that the intersection form of the
resolution $\pi$ is negative definite. More precisely, it is
equivalent to the non-degeneracy of this intersection 
form, as can be easily seen using some diagram
chasing in the cohomology long exact sequence of the pair
$(W,N)$, Poincar{\'e}-Lefschetz duality for the manifold-with-boundary $W$, 
and the fact that $W$ retracts by deformation on $E$. Therefore,
Goresky \& MacPherson's  
theorem is a generalization of the non-degeneracy of the intersection
form associated to a resolution of a normal surface singularity.

\medskip
\section{Constraints on the cohomology of the 
   boundary \\ from the homotopical  
  dimension of the total space} \label{constr}
  
  In this section, all (co)homology groups are considered with coefficients in 
  any commutative ring.  If a space has the homotopical type of a CW-complex 
  of dimension $\leq h$, we say that \emph{it is homotopically of dimension} 
  $\leq h$. Whenever we will be using Poincar{\'e} or Poincar{\'e}-Lefschetz 
  duality morphisms, we will suppose that the orientable manifolds under 
  consideration were arbitrarily oriented. 
  
  \begin{theorem} \label{hyphom}
     Let $W$ be a compact, connected, orientable manifold-with-boundary 
     of dimension $m \geq 4$. Denote by $N$ its boundary. Suppose that 
     $W$ is homotopically of dimension $\leq h$. Consider numbers 
     $i_1,...,i_k \in \{1,..., m-2-h\}$ such that $i_1 +\cdots + i_k
     \geq h +1$.  
     Then the morphism $H^{i_1}(N)\otimes \cdots \otimes H^{i_k}(N) 
     \longrightarrow H^{i_1 +\cdots + i_k}(N)$ induced by the cup-product 
     in cohomology with arbitrary coefficients vanishes identically. 
  \end{theorem}
  
  \textbf{Proof.} By the long exact cohomology sequence of the pair 
     $(W,N)$, we have the exact sequences:
      $$H^{i_l}(W)\stackrel{b^*}{\longrightarrow} H^{i_l}(N) 
            \longrightarrow H^{i_l+1}(W,N), \: \forall \: l \in
            \{1,...,k\},$$ 
    in which $b^*$ is the morphism induced in cohomology by the
    inclusion $N \stackrel{b}{\hookrightarrow} W$. 
    By Poincar{\'e}-Lefschetz duality applied to the oriented 
    manifold-with-boundary $W$, we get $H^{i_l+1}(W,N)\simeq 
    H_{m-i_l-1}(W)$. As $i_l\leq m-2-h$, we see that $m-i_l-1 \geq 
    h+1$. But $W$ was supposed to be homotopically of dimension 
    $\leq h$, therefore $H_{m-i_l-1}(W)=0$. We deduce that all the maps 
    $H^{i_l}(W)\stackrel{b^*}{\longrightarrow} H^{i_l}(N)$ are 
    \emph{surjective}. Consider then the following commutative 
    diagram:
    $$
  \xymatrix{ H^{i_1}(W)\otimes \cdots \otimes H^{i_k}(W) 
         \ar[d]   
             \ar[rr]^{b^*} & 
       & H^{i_1}(N)\otimes \cdots \otimes H^{i_k}(N) 
                 \ar[d]  \\
     H^{i_1 +\cdots + i_k}(W) \ar[rr]^{b^*} &  & H^{i_1 +\cdots + i_k}(N) }$$
As $i_1 +\cdots + i_k \geq h +1$, we get $H^{i_1 +\cdots + i_k}(W) =0$. 
Therefore the composed morphism defined by the commutative diagram 
vanishes identically. But the upper horizontal morphism is surjective, as 
a tensor product of surjective morphisms, therefore the right-hand vertical 
morphism does also vanish identically. \hfill $\Box$

 \begin{remark} \label{condcodim}
    One has to suppose that  $h \leq m-3$ in order to make the set of
    $k$-uples   
    $(i_1,...,i_k)$ satisfying the requested inequalities non-empty
    for some $k\geq 1$. 
  \end{remark}

The next corollary gives restrictions on the cohomology rings of 
Stein fillable manifolds.  In the next section we apply it to give 
examples of Milnor fillable manifolds which are not Stein fillable. 
  
\begin{corollary} \label{Stein}
   Let $N$ be a Stein fillable manifold of dimension $2n-1\geq 5$. 
    Consider numbers 
     $i_1,...,i_k \in \{1,..., n-2\}$ such that $i_1 +\cdots + i_k
     \geq n +1$.    
     Then the morphism $H^{i_1}(N)\otimes \cdots \otimes H^{i_k}(N) 
     \longrightarrow H^{i_1 +\cdots + i_k}(N)$ induced by the cup-product 
     in cohomology with arbitrary coefficients vanishes identically.
   \end{corollary} 

\textbf{Proof.} Combine Theorem \ref{homtype} and Theorem \ref{hyphom}. 
\hfill $\Box$
\medskip

 \begin{remark} 1) We have asked that $N$ be of dimension at least 5, because 
  the previous theorem says nothing about $3$-dimensional Stein fillable 
  manifolds (see Remark \ref{condcodim}). 
 
 2) It would be interesting to find topological properties of 
   Stein fillable manifolds which use in an essential way their 
   \emph{orientation}. The properties stated before are expressed 
   purely in terms of the cohomology
   ring of the manifold, they are not of this type. 

 3) It would also be interesting to find manifolds which admit 
 both a Stein fillable  
 contact structure and a holomorphically fillable but not Stein fillable 
 contact structure. 
  \end{remark}

The next proposition applies Theorem \ref{hyphom} to total spaces of
circle bundles. 
  
\begin{proposition} \label{bundles}
   Let $N\stackrel{p}{\longrightarrow}\Sigma$ be an oriented  
   circle bundle over an 
   orientable closed connected manifold $\Sigma$ of dimension $m-2\geq 2$. 
   Denote by $e \in H^2(\Sigma)$ its Euler class. Consider numbers 
     $i_1,...,i_k \in \{1,..., m-2-h\}$ such that $i_1 +\cdots + i_k
     \geq h +1$.  
     If the morphism $H^{i_1}(\Sigma)\otimes \cdots \otimes H^{i_k}(\Sigma) 
     \longrightarrow H^{i_1 +\cdots + i_k}(\Sigma)$ induced by the cup-product 
     is surjective, but $H^{i_1+\cdots + i_r -2}(\Sigma) 
     \stackrel{\cup e}{\longrightarrow} H^{i_1 +\cdots + i_k}(\Sigma)$ 
     is not surjective, then $N$ does not bound a manifold which is 
     homotopically of dimension $\leq h$. 
 \end{proposition}
    
  \textbf{Proof.} Consider the following part of the Gysin long exact sequence 
    associated to the circle bundle $p$:
    $$ H^{i_1+\cdots + i_r -2}(\Sigma) 
     \stackrel{\cup e}{\longrightarrow} H^{i_1 +\cdots + i_k}(\Sigma) 
      \stackrel{p^*}{\longrightarrow} H^{i_1 +\cdots + i_k}(N)$$ 
      By hypothesis, the morphism on the left is not surjective,
      therefore the one  
      on the right is not vanishing identically.
      
      Look then at the following commutative diagram:
    $$
  \xymatrix{ H^{i_1}(\Sigma)\otimes \cdots \otimes H^{i_k}(\Sigma) 
         \ar[d]   
             \ar[rr]^{p^*} & 
       & H^{i_1}(N)\otimes \cdots \otimes H^{i_k}(N) 
                 \ar[d]  \\
     H^{i_1 +\cdots + i_k}(\Sigma) \ar[rr]^{p^*} &  & H^{i_1 +\cdots +
       i_k}(N) }$$ 
    The left-hand vertical morphism being surjective by hypothesis and the 
    horizontal arrow on the bottom being non-zero, as seen before, 
    we deduce that the morphism defined by the diagram is non-zero. Therefore, 
    the right-hand vertical arrow is non-zero. Theorem \ref{hyphom} allows 
    then to conclude. \hfill $\Box$

\medskip
\section{Applications to Milnor fillable manifolds} 
   \label{appl}

In this section we apply Theorem \ref{hyphom} and Proposition 
\ref{bundles} in order to give lower bounds on the dimensions 
of the exceptional sets of resolutions of isolated singularities 
(Proposition \ref{dimex} and Corollary \ref{corline}), and to construct 
classes of Milnor fillable but not Stein fillable manifolds in any odd 
dimension $\geq 5$ (Proposition \ref{nonsm} and Corollary \ref{corfin}). 
Therefore, we keep working with arbitrary coefficient rings in cohomology. 
\medskip

In Section \ref{notions} we have recalled the notion of divisorial resolution 
of an isolated singularity. We recall now the simplest example of 
normal singularity 
which admits non-divisorial resolutions. Consider the affine quadratic cone
$X$ in $\C^4$,  defined by the equation 
 $xy=zt$. It is the cone over a smooth projective quadric $\Sigma
 \stackrel{s}{\hookrightarrow} \C\gras{P}^3$, 
 which can also be seen as the image of the Segre embedding of 
 $\C\gras{P}^1\times\C\gras{P}^1$, which shows that $\Sigma$ is doubly 
 ruled. The variety $X$ has an isolated singular point at $0$, which  
 can be resolved by blowing it up. The exceptional set is then isomorphic 
 to the projectified tangent cone, that is,  to $\Sigma$:  we have a
 divisorial  resolution.  The 
 total space of this resolution is  isomorphic to the total space
 of the line bundle $s^*\mathcal{O}(-1)$ over $\Sigma$.  
 The initial affine cone $X$  can be seen as obtained by contracting the zero 
 section in this total space. But instead of contracting all of $\Sigma$, 
 one can contract to a point each line in one of the rulings of the quadric, 
 obtaining another resolution whose exceptional set is a smooth 
 rational curve which parametrizes the lines having been contracted. 
 As there are two rulings, one gets two 
 resolutions  with exceptional sets of dimension 1.  The induced 
 birational map between the total spaces of these resolutions is the 
 simplest example of what algebraic geometers call  \emph{ a flop}.  
 
 The following proposition gives topological obstructions to the existence 
 of resolutions with exceptional sets of small dimension. 

\begin{proposition} \label{dimex}
   Let $(X,x)$ be an isolated singularity of normal complex analytic 
   space of complex dimension $n$. Denote by $N$ its abstract boundary. 
   If one can find numbers $i_1,...,i_k \in \{1,..., 2n-2-h\}$ such that 
   $i_1 +\cdots + i_k \geq h +1$ and the morphism 
   $H^{i_1}(N)\otimes \cdots \otimes H^{i_k}(N) 
     \longrightarrow H^{i_1 +\cdots + i_k}(N)$ induced by the cup-product 
     does \emph{not} vanish identically, then the exceptional set of 
     any resolution of $(X,x)$ has complex dimension at least 
     $\frac{h+1}{2}$.
 \end{proposition}
 
 \textbf{Proof.} Consider a resolution $(\tilde{X}, E)\longrightarrow (X,x)$ 
   of the germ $(X,x)$. One can choose
   as representative of $\tilde{X}$ the preimage $W$ 
   of a Milnor neighborhood of $x$ in $X$. Therefore, it retracts by 
   deformation on $E$, which shows that it is homotopically of dimension $2d$, 
   where $d$ is the complex dimension of $E$. By our choice, $W$ 
   has a boundary diffeomorphic to $N$. Theorem \ref{hyphom} applied 
   to our hypotheses 
   implies  that $W$ is not homotopically of dimension $\leq h$, therefore
   $2d \geq h+1$.  \hfill $\Box$ 
   \medskip
 
If $L$ is a line bundle on a projective manifold $\Sigma$, denote by
$N(\Sigma, L)$ the total space of the associated circle bundle.  
By a theorem of Grauert \cite{G 62}, if $L$ is \emph{anti-ample} (that
is, if its dual is ample), one can contract the zero-section
in the total space of $L$, getting like this a normal complex 
affine variety with an isolated singularity as the image of the
zero-section. It is the 
simplest example of Remmert reduction of a holomorphically convex
space (see Peternell \cite[Section 2]{P 94}).  In this case, $N(\Sigma, L)$ is
isomorphic to the boundary of the singularity. The next
corollary applies Proposition \ref{dimex} to the singularities 
obtained in this way.

 \begin{corollary} \label{corline}
   Let $(X,x)$ be an isolated normal singularity obtained by 
   contracting the zero section 
   in the total space of an anti-ample line bundle $L$ over a  
   projective manifold $\Sigma$ of complex dimension $n-1 \geq
   2$. Suppose that   there are numbers  
     $h \in \{2,...,2n-3\}, \: i_1,...,i_k \in \{1,..., m-2-h\}$ such
     that $i_1 +\cdots + i_k   \geq h +1$,  
     the morphism $H^{i_1}(\Sigma)\otimes \cdots \otimes H^{i_k}(\Sigma) 
     \longrightarrow H^{i_1 +\cdots + i_k}(\Sigma)$ induced by the cup-product 
     is surjective, but $H^{i_1+\cdots + i_r -2}(\Sigma) 
     \stackrel{\cup c_1(L)}{\longrightarrow} H^{i_1 +\cdots + i_k}(\Sigma)$ 
     is not surjective. Then the exceptional set of 
     any resolution of $(X,x)$ has complex dimension at least 
     $\frac{h+1}{2}$.
 \end{corollary}

 \textbf{Proof.} This is an immediate consequence of the previous
 proposition and of Proposition \ref{bundles}.   We use the fact that
   the Euler class of the circle bundle $N(\Sigma,L)$ 
   is equal to the first Chern class of $L$. \hfill $\Box$ 
   \medskip

 We have seen before that the contraction of the zero-section of the line  
 bundle $s^*\mathcal{O}(-1)$, where $s$ is the Segre embedding of 
 $\C\gras{P}^1\times\C\gras{P}^1$, admits small resolutions. 
 Consider instead the isolated normal singularity  obtained by contracting the 
 zero-section in the total spaces of the line bundle 
 $s^*\mathcal{O}(-a)$, where  
 $a\geq 2$. Corollary \ref{corline} applied to $L:=s^*\mathcal{O}(-a),
 n=3, h=k=i_1=i_2=2$, implies that this singularity 
 does not admit small resolutions. Indeed, as one can see easily
 using K{\"u}nneth formulae for the product manifold $\Sigma\simeq
 \C\gras{P}^1\times\C\gras{P}^1$, the morphism   
 $H^2(\Sigma)\otimes H^2(\Sigma) 
 \longrightarrow H^4(\Sigma)$ is surjective, but
 the morphism   
 $H^2(\Sigma)\stackrel{\cup
   c_1(L)}{\longrightarrow} H^4(\Sigma)$ is  
 not surjective \emph{if we work with integral cohomology}, as its 
image is divisible by $a \geq 2$. This 
 implies that the exceptional set of any resolution is  
 of complex dimension 2.  
 \medskip

\begin{remark} \label{intcoef}
   The previous example shows that it is essential to apply Corollary
   \ref{corline} using cohomology groups with \emph{integral} coefficients, 
   as the multiplication by the Chern class of
   the line bundle becomes surjective if one uses instead rational or real
   coefficients.  
\end{remark}

\emph{A smoothing} of an isolated singularity is a germ of (flat)
deformation over an irreducible base, and whose generic fiber is smooth. 
If a given singularity admits smoothings, then one says
that \emph{it is smoothable}. 

Consider a smoothing of an isolated normal singularity $(X,x)$. 
By working inside a Milnor representative of the total space of the 
smoothing and using Ehresmann's stability theorem, it can be shown 
that the diffeomorphism type of a generic
fiber is well-defined. It is called \emph{the Milnor fiber} of
the smoothing. Using again Ehresmann's theorem, we see that its boundary 
is diffeomorphic to the boundary of $(X,x)$. 
Recall from Section \ref{notions} that a Milnor representative of the
total space  
of the smoothing is defined as a  sublevel of a (automatically spsh) 
euclidean rug function. Theorem \ref{propstein} implies that this Milnor 
representative is a Stein space. By Theorem \ref{fin}, we see that the
fibers of  
the smoothing are also Stein, therefore the associated Milnor fiber can 
be endowed with the structure of a compact Stein manifold. This 
Stein structure is unique only up to deformations, but this is enough 
in order to see that \emph{the boundary of a smoothable normal isolated  
singularity is Stein fillable}.

Not all isolated singularities are smoothable. 
The following proposition gives smoothing obstructions from the  
cohomology of the boundary. They are different from the ones of   
Hartshorne \cite{H 74}, Rees \& Thomas \cite{RT 78}, Sommese
\cite{S 79} and Looijenga \cite{L 86}. For details on the study of
non-smoothable singularities, 
one can consult Greuel \& Steenbrink  \cite{GS 83}.

\begin{proposition} \label{nonsm}
    Let $(X,x)$ be an isolated singularity of normal complex analytic 
   space of complex dimension $n$. Denote by $N$ its abstract boundary. 
   If one can find numbers $i_1,...,i_k \in \{1,..., n-2\}$ such that 
   $i_1 +\cdots + i_k \geq n +1$ and the morphism 
   $H^{i_1}(N)\otimes \cdots \otimes H^{i_k}(N) 
     \longrightarrow H^{i_1 +\cdots + i_k}(N)$ induced by the cup-product 
     does  \emph{not} vanish identically, then $N$ is not Stein
     fillable. In particular, $(X,x)$ is not smoothable.
 \end{proposition}
 
 \textbf{Proof.} If $(X,x)$ was smoothable then, as explained before,  
 its boundary would be Stein fillable. We get a contradiction from 
 Corollary \ref{Stein}.  \hfill $\Box$
   \medskip
 
 The  following corollaries are examples of application of the previous 
 proposition. Their proofs  use \emph{integral} cohomology.

\begin{corollary} \label{corfin}
    Let $\Sigma$ be an abelian variety of complex dimension $n-1 \geq 2$ 
    and $L$ be an anti-ample line bundle on 
    $\Sigma$ such that $c_1(L)$ is not a primitive element of the lattice 
    $H^2(\Sigma;\Z)$. Then the manifold $N(\Sigma,L)$ is not Stein fillable. 
    In particular, the germ with isolated 
    singularity obtained by contracting the $0$-section of the total
    space of $L$ is not  smoothable. 
 \end{corollary}
 
 \textbf{Proof.}  Topologically, $\Sigma$ is a $(2n-2)$-dimensional
  torus. Therefore, its integral  
  cohomology ring is isomorphic to the exterior algebra $\bigwedge^{2n-2}
  H^1(\Sigma)$. This implies that the morphism 
  $\otimes^{n+1}H^{1}(\Sigma)
     \longrightarrow H^{n+1}(\Sigma)$ is surjective. 
  
  As $c:= c_1(L)$ is not a primitive element of $H^2(\Sigma)$, we can write 
  $c=a\cdot k$,  
  where $a \in \Z, a\geq 2$, and $k\in H^2(\Sigma)$. This implies that 
  $\mbox{im}(H^{n-1}(\Sigma) \stackrel{\cup c}{\longrightarrow}
  H^{n+1}(\Sigma)) \subset a\cdot H^{n+1}(\Sigma)$, which shows that
  the last map  is not surjective.  
  
  Using Theorem \ref{bundles}, we see that all the hypothesis of
  Proposition  \ref{nonsm} are 
  satisfied, with $k=n+1$, $i_1=\cdots=i_k=1$. The
  conclusion  follows. 
\hfill $\Box$

\begin{remark} \label{biran}
    The previous corollary answers partially the following question 
    of Biran \cite{B 05}: ``\emph{Non-fillability of circle bundles $P$ 
    over $\Sigma$ with  $\dim_{\R}\Sigma \geq 4$ would be a new 
    ``contact phenomenon''.  An interesting example to consider seems 
    to be $P \rightarrow \Sigma$, where $\Sigma$ is an Abelian variety 
    of complex dimension $\geq 2$}''. Our answer is partial because we 
    have to impose an hypothesis on $c_1(L)$. 
\end{remark}

 \begin{corollary} \label{HJ}
   Whenever $n\geq 3$ and $a\geq 2$, the manifold $N(\C\gras{P}^{n-1},
   \mathcal{O}(-a))$  
   is not Stein fillable. In particular, the germ with isolated
   singularity obtained by  
   contracting the zero-section of the total space of
   $\mathcal{O}(-a)$ is not smoothable.  
 \end{corollary}

 \textbf{Proof.}   The integral cohomology ring 
   of  $\Sigma:=\C\gras{P}^{n-1}$   is isomorphic to the graded algebra 
   $\Z[x]/(x^n)$, where $\mbox{deg}(x)=2$. This implies that the
   morphism  $\otimes^{n-1} H^2(\Sigma) \longrightarrow 
     H^{2n-2}(\Sigma)$ induced by cup-product, is surjective.

   Moreover, if $c:=c_1(\mathcal{O}(-a))=-ax$, then 
   $\mbox{im}(H^{2n-4}(\Sigma) \stackrel{\cup c}{\longrightarrow} 
  H^{2n-2}(\Sigma)) \subset a\cdot H^{2n-2}(\Sigma)$, which shows that the
  last map  is not surjective, as $  H^{2n-2}(\Sigma)$ is a free group of rank
  one. 
        
  Using Theorem \ref{bundles}, we see that all the hypothesis of
  Proposition  \ref{nonsm} are 
  satisfied, with $k=n-1$, $i_1=\cdots=i_k=2$. The
  conclusion  follows. 
\hfill $\Box$ 
 \medskip

\begin{remark}
  1) The germ obtained by contracting the zero-section of
  $\mathcal{O}(-a)$ has  
  an alternative description as the quotient of $\C^n$ by the
  action of the  
  group of $a$-th roots of unity by coordinate-wise multiplication. Therefore 
  the boundaries of the associated singularities are particular higher
  dimensional lens spaces. 
  For $a=2$, we obtain the real projective space $\R\gras{P}^{2n-1}$. We get 
  like this an alternative proof of the fact that  for $n\geq 3$, this
  space is not Stein fillable (see Eliashberg, Kim \&
   Polterovich \cite[page 1728]{EKP 06}). 

 2) As explained in the previous remark, the germs considered in the
 corollary are particular quotient singularities of dimension at least
 3. Therefore, by a general theorem of Schlessinger \cite{S 71}, they
 are \emph{rigid}, that is, they admit no non-trivial deformations at all. In
 particular they are non-smoothable, which gives an alternative
 proof of the second sentence of the corollary. 
\end{remark}

\medskip
\section{The rational cohomology of Milnor fillable 
   manifolds.} \label{Milnor}

In this section, all cohomology groups are considered \emph{with
  rational coefficients.} This contrasts with the previous section, in
which it was essential to work with integral coefficients (see Remark
\ref{intcoef}).

The next theorem, proved first by Durfee \& Hain \cite{DH 88}, 
after having been announced in \cite{D 86},  states a property  
of the rational cohomology rings  of Milnor fillable manifolds: 

\begin{theorem} \label{obsmiln}
   Let $N$ be a $(2n-1)$-dimensional Milnor fillable manifold, 
  where $n\geq 2$. Consider numbers $\: i_1,...,i_k \in \{1,..., n-1\}$
  such that $i_1 +\cdots + i_k  \geq n$.  
     Then the morphism $H^{i_1}(N)\otimes \cdots \otimes H^{i_k}(N) 
     \longrightarrow H^{i_1 +\cdots + i_k}(N)$ induced by the cup-product 
     in cohomology with rational coefficients vanishes identically. 
\end{theorem}

\textbf{Proof.} Suppose that $N$ is diffeomorphic to the abstract 
boundary of a normal isolated singularity $(X,x)$ of dimension
$n$. Consider a divisorial resolution $(\tilde{X}, E)\longrightarrow (X,x)$ 
of $(X,x)$. Choose as representative of $\tilde{X}$ the preimage $W$ 
of a Milnor representative of $(X,x)$.

Consider  the following commutative diagram, in which the vertical
arrows are induced by the cup-product:

$$
  \xymatrix{ H^{i_1}(W)\otimes \cdots \otimes H^{i_k}(W) 
         \ar[d]   
             \ar[rr]^{b^*} &
       & H^{i_1}(N)\otimes \cdots \otimes H^{i_k}(N) 
                 \ar[d]  \\
     H^{i_1 +\cdots + i_k}(W) \ar[rr]^{b^*} & & H^{i_1 +\cdots +
       i_k}(N) }$$ 

Using point 4) of Theorem     \ref{decomp}, the hypotheses made on the numbers 
$i_1,...,i_k$ imply that the upper
horizontal morphism is surjective as a tensor product of surjective morphisms 
and that the lower horizontal one  
vanishes identically. The conclusion follows. \hfill $\Box$
\medskip

Note that, although the theorems \ref{hyphom} and \ref{obsmiln} are formally
 similar, the key results used in their proofs are completely different. 
 
 In fact, Durfee \& Hain stated the previous theorem for $k=2$. The
 precise bounds  
 $n-1$ and $n$ for $i_1,...,i_k$ and $i_1 +\cdots +i_k$ respectively,
 make the two  
 versions equivalent. They proved their theorem using directly the theory of 
 Beilinson, Bernstein, Deligne \& Gabber, that is, without passing through the 
 theorem of Goresky \& MacPherson. They get also restrictions on the 
 boundaries of tubular neighborhoods of higher dimensional subspaces of 
 algebraic varieties. 

\begin{remark} \label{remiln}
  1) By using Poincar{\'e} duality, the previous theorem may be reformulated 
  in the following way: \emph{on a $(2n-1)$-dimensional Milnor
    fillable manifold,  
  the intersection number of rational homology classes of dimension 
  at least   $n$ is equal to zero}. Of course, we suppose that we take 
  classes whose sum of codimensions is equal to $2n-1$, in order to have 
  a well-defined intersection number.

  2) The previous theorem was obtained in the particular case $n=2, k=2,
  i_1=i_2=1$ by Sullivan \cite{S 75}. In this case, where $(X,x)$ is a
  germ of surface, it results immediately from the fact that the
  intersection form of a resolution of the isolated singular point
  under study is non-degenerate. As explained in Section
  \ref{notions}. Theorem \ref{decomp} is an
  analog of this result in higher dimensions. From his theorem, 
  Sullivan deduced that the
  boundary of an isolated singularity cannot be diffeomorphic to a
  3-dimensional torus. The following generalizes this to all
  dimensions. 
\end{remark}

As an immediate application of the previous theorem, we see that 
 \emph{for all $n\geq 2$, the torus $\gras{T}^{2n-1}$ is not 
 Milnor-fillable}. Indeed, as the cohomology ring of $\gras{T}^{2n-1}$ 
is isomorphic to the exterior  algebra $\bigwedge^{\bullet} 
  H^1(\gras{T}^{2n-1})$ (a fact already used in the proof of Corollary
  \ref{corfin}),  the morphism 
  $\otimes^{n}H^{1}(\gras{T}^{2n-1}) \longrightarrow H^{n}(\gras{T}^{2n-1})$ is
  surjective and does not vanish identically. We apply then 
  the previous theorem to $k=n, i_1=\cdots=i_n=1$. In the next section 
    we will see that  starting from dimension 5, odd-dimensional 
    tori are not even holomorphically fillable. 

 The fact that for $n\geq 2$, the torus $\gras{T}^{2n-1}$ is not the 
        boundary of an \emph{isolated complete intersection singularity}
        (abbreviated \emph{icis})  
        can be proved differently. Indeed, it is implied by the fact that 
        the boundary of an icis of complex dimension $n$ is
        $(n-2)$-connected.  
        This was proved by Milnor \cite{M 68} for hypersurfaces and 
        generalized by Hamm \cite{H 71} to arbitrary icis. 
        
      As noted by Durfee \cite{D 86}, Theorem \ref{obsmiln} implies more 
     generally that no manifold of the form $K_1 \times K_2 \times K_3$ with 
     $\mbox{dim} \: K_1 + \mbox{dim} \: K_2 + \mbox{dim} \: K_3 = 2n-1$ 
     and $\mbox{dim} \: K_i \leq n-1, \: \forall \: i \in \{1,2,3\}$
     is Milnor fillable.

\medskip
\section{The rational cohomology of holomorphically 
 fillable manifolds.} 
\label{holomorph}

In what follows, all (co)homology groups are considered with \emph{rational} 
coefficients. 
  
  \begin{proposition} \label{surgery}
     Suppose that $n \geq 3$. Let $W$ be an oriented cobordism of
     dimension $2n$  
     from a manifold $N_1$ to a manifold $N_2$, such that $(W, N_1)$ 
     has the homotopy type of a 
     relative CW-complex of dimension $\leq n$. Consider numbers 
     $i_1,...,i_k \in \{1,..., n-2\}$ such that $i_1 +\cdots + i_k \geq n +1$. 
     If the morphism $H^{i_1}(N_1)\otimes \cdots \otimes H^{i_k}(N_1) 
     \longrightarrow H^{i_1 +\cdots + i_k}(N_1)$ induced by the cup-product 
     in cohomology with rational coefficients vanishes identically,
     then the same is true for the analogous  
     morphism associated to $N_2$. 
  \end{proposition}
  
  \textbf{Proof.} For $j=1,2$, denote
  by $N_j \stackrel{u_j}{\hookrightarrow}W$ the 
  inclusion morphism.  Our hypothesis   on the pair $(W, N_1)$ implies that: 
  $$H^i(W,N_1)\simeq 0, \: \forall \: i \geq n+1.$$

  Using the exact cohomology sequence of the pair $(W, N_1)$, we deduce that:
  \begin{equation} \label{isom}
    H^i(W)\stackrel{u_1^*}{\longrightarrow}H^i(N_1) \mbox{ is injective}, 
    \: \forall \: i \in \{n+1,..., 2n-1\},
  \end{equation}
 (in fact those morphisms are bijective, but for the  proof of the proposition 
 we need only their injectivity).
 
  By generalized Poincar{\'e}-Lefschetz duality applied to the cobordism
  $W$ from $N_1$  
  to $N_2$ (see Dold \cite[page 307]{D 80}), we get $H^{i+1}(W, N_2) 
  \simeq$ \linebreak $H_{2n-i-1}(W, N_1)$, therefore: 
  $$H^{i+1}(W, N_1)\simeq 0, \: \forall \: i \in \{0,...,n-2\}.$$
 Using the exact cohomology sequence of the pair $(W, N_2)$, we 
  deduce that:
   \begin{equation} \label{surj}
    H^i(W)\stackrel{u_2^*}{\longrightarrow}H^i(N_2) \mbox{ is surjective},
    \: \forall \: i \in    \{0,...,n-2\}.
  \end{equation}
  
  Consider then the following commutative 
    diagram:
    \begin{scriptsize} 
     $$
  \xymatrix{ H^{i_1}(N_1)\otimes \cdots \otimes H^{i_k}(N_1) 
                 \ar[d]   & \ar[l]_{u_1^*} 
    H^{i_1}(W)\otimes \cdots \otimes H^{i_k}(W) 
         \ar[d]   
             \ar[r]^{u_2^*} 
       & H^{i_1}(N_2)\otimes \cdots \otimes H^{i_k}(N_2) 
                 \ar[d]  \\
        H^{i_1 +\cdots + i_k}(N_1)   & \ar[l]_{u_1^*} 
     H^{i_1 +\cdots + i_k}(W) \ar[r]^{u_2^*} &   H^{i_1 +\cdots + i_k}(N_2) }$$
     \end{scriptsize}
  By hypothesis, the left-side vertical arrow vanishes. As the left-side 
  lower horizontal arrow is injective by (\ref{isom}), we deduce that 
  the middle arrow vanishes identically. As the right-side upper horizontal 
  arrow is surjective by (\ref{surj}), we deduce that the right-side vertical 
  arrow also vanishes. 
 \hfill $\Box$
 
 \medskip
 As a consequence, we get the following property of the rational 
 cohomology rings of holomorphically fillable manifolds, proved by Bungart 
 \cite{B 92}:

  \begin{theorem} \label{holofill}
     Suppose that $n \geq 3$. Let $N$ be  a holomorphically fillable  manifold 
     of dimension $2n-1$.  Consider numbers 
     $i_1,...,i_k \in \{1,..., n-2\}$ such that $i_1 +\cdots + i_k \geq n +1$. 
     Then the morphism $H^{i_1}(N)\otimes \cdots \otimes H^{i_k}(N) 
     \longrightarrow H^{i_1 +\cdots + i_k}(N)$ induced by the cup-product 
     in cohomology with rational coefficients vanishes identically. 
  \end{theorem}

 \textbf{Proof.} Suppose that $N$ is the strictly pseudo-convex boundary of 
 a compact holomorphic manifold $Z$ of complex dimension $2n$. Denote
 by $Z \stackrel{r}{\rightarrow} X$ the Remmert reduction morphism of
 $Z$. Then $X$ is a Stein space with at most isolated singularities,
 obtained as images of some  of the connected components of the maximal
 compact analytic subset of $Z$ (the other connected components contract to 
 smooth points of $X$). 

Denote by $F$ the (finite) set of
 singular points of $X$.  Let $\phi$ be a function as described in the
 conclusion of Proposition \ref{specfunc}.  Suppose (without reducing
 the generality) that the absolute minimum of $\phi$, attained by 
 hypothesis on the finite set $F$, is equel to $0$ (it is enough
 to add a constant to $\phi$ in order to get this). Suppose
 moreover that $\phi$ is a Morse function with only one critical point on
 each critical level, which can be realized by a sufficiently small smooth
 perturbation on a compact subset of $X$ (we use the stability of 
 spsh functions formulated in Proposition \ref{spsh}, point 3)). For
 $\epsilon>0$ sufficiently small, the level $N:= \phi^{-1}(\epsilon)$ is a
 disjoint union of manifolds diffeomorphic to the boundaries of the
 singularities of $X$.

     Denote by $\nu>0$ the maximum value of $\phi$, attained by construction 
  exactly on $N_2 := N$. 
     Then   $W:=\rho^{-1}([\epsilon, \nu])$ is an oriented cobordism from 
     $N_1$ to $N_2$. As 
  $\phi$ is spsh, we deduce from Proposition \ref{index} that $(W, N_1)$ has 
  the  homotopy type of a relative CW-complex of dimension $\leq n$. By 
  Theorem \ref{obsmiln}, for each connected component $C$ of $N_1$, 
  the morphism $H^{i_1}(C)\otimes \cdots \otimes H^{i_k}(C) 
     \longrightarrow H^{i_1 +\cdots + i_k}(C)$ vanishes identically, which 
     implies that the same is true for $N_1$. Therefore we can apply 
     Proposition \ref{surgery}, and the conclusion follows. 
  \hfill $\Box$

 \medskip
 
 A similar theorem holds for Stein fillable manifolds, with the  
 essential difference that it is then true for cohomology groups 
 \emph{with arbitrary coefficient rings} (Theorem \ref{Stein}). This 
 difference allows to detect with our methods holomorphically fillable 
 manifolds which are not Stein fillable, as we did in Section \ref{appl}.  
  
 \begin{remark} 
 1) By using Poincar{\'e} duality, the previous theorem may be reformulated 
  in the following way: \emph{on a $(2n-1)$-dimensional 
  holomorphically  fillable manifold, 
  the intersection number of rational homology classes of dimension 
  at least   $n+1$ is equal to zero} (compare with Remark \ref{remiln}, 1)).

 2) One has to introduce the restriction $n\geq 3$ in order 
 to have integers  $i_1,..., i_k$ which satisfy the conditions of the
 hypotheses.  
Therefore, as in the case of Corollary \ref{Stein}, the previous theorem says 
nothing about 3-dimensional manifolds. 

  3) With the notations of the proof of Theorem \ref{holofill}, the
  surjectivity of the morphism
  $H^0(W)\stackrel{u_2^*}{\longrightarrow}H^0(N_2)$ (a consequence of
  (\ref{surj})),  
  shows that the boundary of a strongly pseudoconvex
  \emph{connected} manifold is also connected (folklore). 
 \end{remark}
 
 The proof of the previous theorem shows that one can get more
 information on the  
 cohomology rings of holomorphically fillable manifolds from more detailed 
 knowledge of the topology of the isolated singularities of a Stein
 space which fills it.  
 For example:
 
 \begin{proposition} 
    Let $N$ be the boundary of a compact Stein space $X$ with isolated
    singularities.  
    Fix a ring of coefficients $A$ and numbers $i_1,..., i_k \in
    \{1,...,n-2\}$ such that  
    $i_1+\cdots + i_k\geq n+1$. If the morphism $H^{i_1}(M,A)\otimes
    \cdots \otimes    H^{i_k}(M,A)\rightarrow  
    H^{i_1+\cdots + i_k}(M,A)$ induced by the cup-product vanishes
    identically for the abstract    boundary $M$ of  
    each isolated singular point of $X$, then the same is true for $N$. 
 \end{proposition}
 
 As an immediate application of Theorem \ref{holofill}, we see that 
  \emph{for all $n \geq 3$, the torus $\T^{2n-1}$ is not
    holomorphically fillable.}  
 Indeed, one has simply to apply the previous theorem to 
 $k=n+1, i_1=\cdots =i_{n+1}=1$. 
 
Note that the torus $\T^{2n-1}$ can be realised as a Levi-flat boundary of
  a complex  
       manifold: consider the product of an abelian variety and the
       closed unit  
       disc in $\C$. This shows the importance of the \emph{strong} 
       convexity hypothesis in the definition of holomorphically
       fillable manifolds.  
  
  By a theorem of Bourgeois \cite{B 02} (see also Giroux \cite{G
    02} for details on   
  the context of research having led to it), if a closed orientable
  manifold $M$ 
  admits a contact structure, then $M\times \T^2$ also does. This
  implies that   
  all odd-dimensional tori  admit contact structures, as $\T^3$ does 
  (see the next paragraph). The previous corollary shows that a contact 
   structure on a torus of dimension at least $5$ cannot be
   holomorphically fillable.

 Instead, $\T^3$ is holomorphically fillable: it can be realized as a 
 strongly pseudoconvex boundary of a tubular neighborhood of 
 $\gras{S}^1\times \gras{S}^1$ 
 standardly embedded in $\C^2$ (see Eliashberg \cite{E 96}). 
 By the theorem of Sullivan quoted in the previous section and 
 generalized in Theorem \ref{obsmiln}, $\T^3$ is not Milnor 
 fillable. In a similar way, we get using Theorem \ref{obsmiln}:  
 
 \begin{proposition} \label{prod}
   For any $n\geq 2$, the product $\T^n \times \gras{S}^{n-1}$ is
   holomorphically fillable but not Milnor fillable. 
 \end{proposition}
 
    \textbf{Proof.} Consider the standard embedding of $\T^n=
    \gras{S}_1^1\times \cdots  
  \times \gras{S}_n^1$ in $\C^n$ as the product of unit circles in
  each factor $\C$ (the indices denote different copies of $\gras{S}^1$). 
  As the image of this embedding is totally real, we see  
  that it has strongly pseudoconvex regular neighborhoods 
  (see Grauert \cite{G 58},  Eliashberg \cite{E 97}). The  boundaries
  of these regular   
  neighborhoods are diffeomorphic to $\T^n\times \gras{S}^{n-1}$,
  which shows that this  
  last manifold is holomorphically fillable. 
  
  Choose now points $P_i \in \gras{S}_i^1, \: \forall \: i \in
  \{1,...,n\}$ and  
  $P \in \gras{S}^{n-1}$. The submanifolds $K_1:=\gras{S}_1^1\times P_2 \times 
  \cdots \times P_n \times \gras{S}^{n-1}, \: K_2:= \T^n \times P, 
  \: K_3:= P_1 \times \gras{S}_2^1\times \cdots \times
  \gras{S}_n^1\times \gras{S}^{n-1}$ of  
  $\T^n\times \gras{S}^{n-1}$ have only the point $P_1\times\cdots
  \times P_n \times P$    in common, where they meet 
  transversally. Therefore, with convenient choices of orientation,
  the intersection number  
  of their homology classes is equal to $1$. For $j \in \{1,2\}$, denote by 
  $\gamma_j \in H^{n-1}( \T^n\times\gras{S}^{n-1})$ the Poincar{\'e} dual of the 
  homology class of $K_j$. We deduce that $\gamma_1\cup \gamma_2$ does
  not vanish  
  in $H^{2n-2}(\T^n\times\gras{S}^{n-1})$, which shows that the morphism 
  $H^{n-1}( \T^n\times\gras{S}^{n-1})\otimes H^{n-1}( \T^n\times\gras{S}^{n-1})
  \longrightarrow H^{2n-2}( \T^n\times\gras{S}^{n-1})$ induced by the
  cup-product  
  does not vanish identically. By Theorem \ref{obsmiln}, we deduce that 
  $\T^n\times \gras{S}^{n-1}$ is not Milnor fillable. \hfill $\Box$

\medskip


{\small


\begin{thebibliography}{99}

\bibitem{AF 59} Andreotti, A., Frankel, T. \textit{The Lefschetz theorem 
  on hyperplane sections.} Annals of Maths. \textbf{69}, no.3 (1959), 
  713-717.

\bibitem{BBD 82} Beilinson, A.A., Bernstein, J. N., Deligne,
  P. \textit{Faisceaux pervers.} Ast{\'e}risque \textbf{100},
  Soc. Math. France, 1982. 
  
\bibitem{BN 90} Bell, S.R., Narasimhan, R. \textit{Proper Holomorphic 
  Mappings of Complex Spaces.} Several Complex Variables VI, 
  Encyclopaedia of Maths. Sciences vol. \textbf{69}, Barth, W.,
  Narasimhan, R. eds.,  Springer Verlag 1990.

\bibitem{B 90} Bennequin, D. \textit{Topologie symplectique, convexit{\'e} 
  holomorphe et structures de contact [d'apr{\`e}s Y. Eliashberg, D. 
  Mc Duff et al.]} S{\'e}minaire Bourbaki no. \textbf{725}, Ast{\'e}risque 
  189-190, Soc. Math. France, 1990. 

\bibitem{B 05} Biran, P. \textit{Symplectic topology and algebraic families.} 
   Proc. of the European Congress of Maths. (Stockholm 2004), European 
   Math. Soc., 2005. 

\bibitem{BO 97} Bogomolov, F.A., de Oliveira, B. \textit{Stein Small
    Deformations of Strictly Pseudoconvex Surfaces.} Contemporary
  Mathematics \textbf{207} (1997), 25-41. 
  
  \bibitem{B 59} Bott, R. \textit{On a theorem of Lefschetz.} Michigan Math. 
  Journal \textbf{6} (1959), 211-216. 

\bibitem{B 02} Bourgeois, F. \textit{Odd dimensional tori are contact
  manifolds.}  Int. Math. Res. Not. \textbf{30} (2002), 1571--1574.  

\bibitem{B 92} Bungart, L. \textit{Vanishing cup products on
    pseudoconvex CR manifolds.} Contemp. Math. \textbf{137} (1992),
  105-111. 

\bibitem{CM 05} de Cataldo, M.A.A., Migliorini, L. \textit{The Hodge
    theory of algebraic maps.} Ann. Sci. ENS \textbf{38} (2005),
  693-750. 
  

\bibitem{CNP 06} Caubel, C., N{\'e}methi, A., Popescu-Pampu,
  P. \textit{Milnor open books and Milnor fillable contact
    3-manifolds.}   Topology \textbf{45} (2006), 673-689. 

\bibitem{D 80} Dold, A. \textit{Lectures on Algebraic Topology.}
  Springer, 1980. 

\bibitem{D 86} Durfee, A.H. \textit{Topological restrictions on the
    links of isolated complex singularities.} In \textit{Complex
    analytic singularities.} Adv. Studies in Pure Maths. \textbf{8}
  (1986), 95-99. 

\bibitem{DH 88} Durfee, A.H., Hain, R.M. \textit{Mixed Hodge
    structures on the homotopy of links.} Math. Ann. \textbf{280}
  (1988), 69-83.
    
\bibitem{E 89} Eliashberg, Y. \textit{Filling by holomorphic discs and its 
  applications.}   In \textit{Geometry of low-dimensional manifolds. 2.} 
  (Durham, 1989), 45-67, London Math. Soc. Lecture Notes Series 
  \textbf{151}, Cambridge Univ. Press, 1990.

\bibitem{E 90} Eliashberg, Y. \textit{Topological characterization of
    Stein manifolds of dimension >2.} Internat. Journ. of
  Maths. \textbf{1} no. 1 (1990), 29-46. 
  
\bibitem{E 96} Eliashberg, Y. \textit{Unique holomorphically fillable contact 
   structure on the 3-torus.}  Int. Math. Res. Not. \textbf{2}
 (1996), 77-82.  
  
  
\bibitem{E 97} Eliashberg, Y. \textit{Symplectic geometry of 
  plurisubharmonic functions.} (Notes by M. Abreu) In 
   \textit{Gauge theory and symplectic geometry.} J. Hurtubise and 
  F. Lalonde eds., 49-67, Kluwer, 1997. 
    
 \bibitem{EG 91} Eliashberg, Y., Gromov, M. \textit{Convex symplectic 
    manifolds.} In \textit{Several complex variables and complex 
    geometry.} Part 2, Proc. Sympos. 
    Pure Math. \textbf{52}, 135-162 (1991). 

\bibitem{EKP 06} Eliashberg, L., Kim, S.S., Polterovich,
  L. \textit{Geometry of contact transformations and domains:
    orderability versus squeezing.} Geometry \& Topology \textbf{10}
  (2006), 1635-1747. 
    

\bibitem{G 08} Geiges, H. \textit{An introduction to contact
    topology.} Cambridge Univ. Press, 2008.


\bibitem{G 02} Giroux, E. \textit{G{\'e}om{\'e}trie de contact: de la
    dimension trois vers les dimensions sup{\'e}rieures.}  Proceedings of
    the International Congress 
    of Mathematicians, Vol. II (Beijing, 2002),  405-414, Higher
    Ed. Press, Beijing, 2002. 
    
\bibitem{GM 82} Goresky, M., MacPherson, R. \textit{On the topology of 
  complex algebraic maps.} In \textit{Algebraic Geometry. Proceedings, 
  La R{\'a}bida, 1981.} J.M. Aroca, R. Buchweitz, M. Giusti, M. Merle eds. 
  L.N.M. \textbf{961}, Springer Verlag, 1982.      
    
\bibitem{G 58} Grauert, H. \textit{On Levi's problem and the imbedding of 
  real-analytic manifolds.} Annals of Maths. \textbf{68} (1958), 460-472.     


\bibitem{G 62} Grauert, H. \textit{{\"U}ber Modifikationen  und exceptionelle
   analytische Mengen.} Math. Ann. {\bf 146} (1962), 331-368.
   
\bibitem{GR 79} Grauert, H., Remmert, R. \textit{Theory of Stein spaces.} 
   Springer-Verlag, 1979. 

   
 \bibitem{GS 83} Greuel, G.-M., Steenbrink, J. \textit{On the topology
     of smoothable  
    singularities.} Proc. of Symp. in Pure Maths. \textbf{40} (1983),
  Part 1, 535-545. 
   
 \bibitem{G 85} Gromov, M. \textit{Pseudoholomorphic curves in 
   symplectic manifolds.} Inv. Math. \textbf{82} (1985), 307-347.

\bibitem{H 71} Hamm, H. \textit{Lokale topologische Eigenschaften
    komplexer R{\"a}ume.} Math. Ann. \textbf{191} (1971), 235-252.

   
 \bibitem{H 86} Hamm, H. \textit{Zum homotopietyp $q$-vollst{\"a}ndiger 
   R{\"a}ume.} J. Reine Angew. Math. \textbf{364} (1986), 1-9. 
   
\bibitem{H 74} Hartshorne, R. \textit{Topological conditions for
    smoothing algebraic singularities.} Topology \textbf{13} 
  (1974), 241-253. 
  
  \bibitem{K 79} Karchyauskas, K.K. \textit{Homotopy properties of complex 
   algebraic sets.} Studies in topology, Steklov Institute (Leningrad), 1979. 


\bibitem{L 86} Looijenga, E. \textit{Riemann-Roch and smoothings of
    singularities.} Topology \textbf{25} (1986), 293-302. 
  

\bibitem{M 63} Milnor, J. \textit{Morse theory.} Princeton Univ. Press, 1963.

\bibitem{M 68} Milnor, J. \textit{Singular Points of Complex
    Hypersurfaces.} Princeton Univ. Press, 1968.
    
\bibitem{M 61} Mumford, D. \textit{The Topology of Normal
    Singularities of an Algebraic Surface and a Criterion for
    Simplicity.} Publ. Math. IHES \textbf{9} (1961), 229-246.



\bibitem{N 67} Narasimhan, R. \textit{On the homology groups of Stein spaces.} 
    Inv. Math. \textbf{2} (1967), 377-385. 

    
\bibitem{N 85} Navarro Aznar, V. \textit{Sur la th{\'e}orie de Hodge des 
   vari{\'e}t{\'e}s alg{\'e}briques {\`a} singularit{\'e}s isol{\'e}es.} In
 \textit{Differential systems and singularities.} (Luminy, 1983),
 Ast{\'e}risque \textbf{130} (1985), 272-307. 

\bibitem{P 94} Peternell, Th.  \textit{Pseudoconvexity, the Levi problem and 
    vanishing theorems.} Several Complex Variables VII,   
     Encyclopaedia of Maths. Sciences vol. 69, Grauert, H., Peternell,
     Th., Remmert, R. eds.,    Springer-Verlag 1994, 221-257.


\bibitem{RT 78} Rees, E., Thomas, E. \textit{Cobordism obstructions to 
   deforming isolated singularities.} Math. Ann. \textbf{232} 
   (1978), 33-54.
      
\bibitem{S 71} Schlessinger, M. \textit{Rigidity of quotient
    singularities.} Inv. Math. \textbf{14} (1971), 17-26. 
    
 \bibitem{S 53} Serre, J.-P. \textit{Quelques probl{\`e}mes globaux relatifs aux 
   vari{\'e}t{\'e}s de Stein.} Colloque sur les fonctions de plusieurs variables, 
   Bruxelles (1953), 57-68.

\bibitem{S 79} Sommese, A.J. \textit{Non-smoothable varieties.} Comment. 
  Math. Helv. \textbf{54} (1979), 140-146.

\bibitem{S 83} Steenbrink, J. \textit{Mixed Hodge structures
    associated with isolated singularities.} Proc. of Symposia in Pure
  Maths. \textbf{40} (1983), Part 2, 513-536. 

\bibitem{S 87} Steenbrink, J. \textit{Mixed Hodge structures and
    singularities: a survey.} In \textit{G{\'e}om{\'e}trie alg{\'e}brique et
    applications. III.} Comptes-rendus de la deuxi{\`e}me conf{\'e}rence
  internationale de La R{\'a}bida (1984), J-M. Aroca,
  T. Sanchez-Giralda, J-L. Vicente eds., Hermann, 1987. 

\bibitem{S 75} Sullivan, D. \textit{On the intersection ring of
    compact three manifolds.} Topology \textbf{14} (1975), 275-277. 
    
\bibitem{T 59} Thom, R. \textit{Les structures diff{\'e}rentiables des boules 
  et des sph{\`e}res.} In \textit{Colloque de g{\'e}om{\'e}trie diff{\'e}rentielle 
  globale.} Bruxelles, dec. 1958. Centre Belge de Recherches 
  Math{\'e}matiques, 1959. 
  
\bibitem{V 44} Du Val, P. \textit{On absolute and non-absolute
    singularities of algebraic surfaces.} Revue de la Facult{\'e} des
    Sciences de l'Univ. d'Istanbul (A) \textbf{91} (1944), 159-215.

 
\bibitem{V 80} Varchenko, A.N. \textit{Contact structures and isolated
    singularities.} Mosc. Univ. Math. Bull. \textbf{35} no.2 (1980),
  18-22.






\end{thebibliography} }
\end{document}